\documentclass{amsart}
\usepackage[dvips]{graphics}
\usepackage{amsfonts}
\usepackage{amssymb}
\usepackage{comment}

\DeclareMathAlphabet{\eusm}{U}{}{}{}  % Euler script math
\SetMathAlphabet\eusm{normal}{U}{eus}{m}{n}
\SetMathAlphabet\eusm{bold}{U}{eus}{b}{n}

\DeclareMathAlphabet{\eufrak}{U}{}{}{}  % Euler fraktur math
\SetMathAlphabet\eufrak{normal}{U}{euf}{m}{n}
\SetMathAlphabet\eufrak{bold}{U}{euf}{b}{n}

\newtheorem{theorem}{Theorem}[section]
\newtheorem{proposition}[theorem]{Proposition}
\newtheorem{lemma}[theorem]{Lemma}
\newtheorem{corollary}[theorem]{Corollary}
\newtheorem{problem}[theorem]{Problem}

\theoremstyle{definition}
\newtheorem{definition}[theorem]{Definition}
\newtheorem{example}[theorem]{Example}

\theoremstyle{remark}
\newtheorem{remark}[theorem]{Remark}

\numberwithin{equation}{section}

\newenvironment{rlist}
{

\begin{enumerate}}
{\end{enumerate}}

\begin{document}

\title{On Square Roots of the Haar State on Compact Quantum Groups}
\author{Uwe Franz}
\address{UF: D\'epartement de math\'ematiques de Besan\c{c}on,
Universit\'e de Franche-Comt\'e 16, route de Gray, 25 030
Besan\c{c}on cedex, France}
\email{uwe.franz@univ-fcomte.fr}
\urladdr{http://www-math.univ-fcomte.fr/pp\_Annu/UFRANZ/}
\thanks{U.F.\ was supported by a Marie Curie Outgoing International
Fellowship of the EU (Contract Q-MALL MOIF-CT-2006-022137) and an ANR Project
(Number 2011 BS01 008 01).}

\author{Adam Skalski}
\address{AS: Institute of Mathematics of the Polish Academy of Sciences,
ul.\'Sniadeckich 8, 00-956 Warszawa, Poland} \email{a.skalski@impan.pl}
\urladdr{http://www.impan.pl/~skalski}

\author{Reiji Tomatsu}
\address{RT: Department of Mathematics
Tokyo University of Science, Yamazaki 2641, Noda, Chiba, 278-8510, Japan}
\email{tomatsu\_reiji@ma.noda.tus.ac.jp}
\urladdr{http://www.ma.noda.tus.ac.jp/u/rto/index.html}

\begin{abstract}
The paper is concerned with the extension of the classical study of probability measures on a compact group which
are square roots of the Haar measure, due to Diaconis and Shahshahani, to the context of compact quantum groups.
We provide a simple characterisation for compact quantum groups which admit no non-trivial square roots of the
Haar state in terms of their corepresentation theory. In particular it is shown that such compact quantum groups
are necessarily of Kac type and  their subalgebras generated by the coefficients of a fixed  two-dimensional
irreducible corepresentation are isomorphic (as finite quantum groups) to the algebra of functions on the group
of unit quaternions. An example of a quantum group whose Haar state admits no nontrivial square root and which is
neither commutative nor cocommutative is given.
\end{abstract}
 \maketitle

\section{Introduction}
\label{sec-intro}

If $G$ is a (locally) compact group, then the space $\textup{M}(G)$ of all bounded regular measures on $G$ is
equipped with a natural product $\star$ afforded by the convolution, making $\textup{M}(G)$ a Banach algebra. In
particular a convolution of two probability measures remains a probability measure. Convolution equations in
$\textup{M}(G)$ have often natural interpretations -- for example the Haar measure $\mu_G$ on compact $G$ can be
described as a unique probability measure in $\textup{M}(G)$ such that $\mu_{G} \star \mu = \mu \star
\mu_G=\mu_G$ for all $\mu \in \textup{Prob}(G)$, idempotents in $\textup{Prob}(G)$ can be characterised as Haar
measures on compact subgroups of $G$ (Kawada-It\^o Theorem), and so on. In  \cite{diaconis+shahshahani86}
Diaconis and Shahshahani showed that the Haar measure of a separable compact  topological group  $G$ does not
admit a non-trivial square root, i.e.\ a probability $\nu\not=\mu_G$ with $\nu\star\nu=\mu_G$, if and only if $G$
is  abelian or of the form $H\times E$, where $H$ is the eight element group of unit quaternions and $E$ a
product of two element groups.

If $(\mathsf{A}, \Delta)$ is a compact quantum group in the sense of Woronowicz (\cite{woronowicz98}), then the
quantum counterpart of the set of the probability measures on the group is given by the state space of
$\mathsf{A}$, $S(\mathsf{A})$. It is again equipped with a natural convolution operation and it makes sense to
ask for the solutions of analogous equations as those listed above. Thus for example the Haar state can be
described as a unique state $h$ on $\mathsf{A}$ such that for all $\rho \in S(\mathsf{A})$ there is $h\star
\rho=\rho \star h = h$, but on the other hand there may exist idempotent states on $\mathsf{A}$ which do not arise
as Haar states on compact quantum subgroups of $\mathsf{A}$ (see \cite{franz+skalski-idemp},
\cite{franz+skalski+tomatsu09}, \cite{Salmi+Skalski2011}).

In this paper we investigate the quantum counterpart of the question studied by Diaconis and Shahshahani -- which
compact quantum groups belong to the \emph{quantum $DS$-family}, i.e.\ have the property that their Haar state
does not admit a non-trivial square root? The proof of the characterisation in \cite{diaconis+shahshahani86}
consists of three main steps: first they show that that the existence of the nontrivial square roots of the Haar
measure on $G$ is equivalent to the existence of a non-zero bounded real nilpotent measure on $G$, then deduce
from that that if $G$ admits no such nilpotent measures then it must be hamiltonian (i.e.\ all its closed
subgroups are normal) and finally classify compact separable hamiltonian groups and use this classification to
complete the proof.

Here we first show that the existence of the nontrivial square roots of the Haar state is equivalent to the fact
that the dual discrete algebraic quantum group $(\hat{\mathcal{A}},\hat{\Delta})$ contains a hermitian non-zero
nilpotent element, cf.\ Theorem \ref{thm-nilp} (the proof in our context is similar to that of
\cite{diaconis+shahshahani86}, but substantially more technical).   Motivated by our work on idempotent states,
we propose a definition of a \emph{hamiltonian compact quantum group} as the one on which all idempotent states
are central and show that if $(\mathsf{A}, \Delta)$ belongs to the quantum $DS$-family, then it is hamiltonian,
cf. Proposition \ref{prop-normal} (in particular, all compact quantum subgroups of $(\mathsf{A},\Delta)$ are
normal).  As the classification of hamiltonian compact quantum groups is currently beyond our reach, to continue
the investigation we need to provide another strategy, based on the corepresentation theory. It turns out that
the quantum $DS$-family consists  exactly of those compact quantum groups which have only one- and
two-dimensional irreducible corepresentations and satisfying a certain additional condition on the linear
functionals coming from their two-dimensional corepresentations, cf.\ Theorem \ref{thm-structure}. This allows us
to deduce that if $(\mathsf{A}, \Delta)$ is in the quantum $DS$-family, then it is necessarily of Kac type and
its subalgebras generated by the coefficients of a fixed  two-dimensional irreducible corepresentation are
isomorphic to the algebra of functions on the group of unit quaternions. That result and some further
observations on the interaction between two-dimensional and one-dimensional corepresentations  are used to
provide an explicit example of a compact quantum group in the quantum $DS$-family which is neither commutative
nor cocommutative.

The detailed plan of the paper is as follows: Section 2 contains all the preliminary facts and terminology
related to compact (and discrete) quantum groups and their corepresentations. Here we also introduce the
fundamental notion of the square root of the Haar state and characterise commutative and cocommutative elements
of the quantum $DS$-family.   Section 3 is devoted to establishing the equivalence between the existence of
nontrivial square roots and non-zero bounded hermitian nilpotent functionals and discus hamiltonian compact
quantum groups. In Section 4 the corepresentation theory starts to play a prominent role, providing a means to
characterise the quantum $DS$-family. This is used in the following section to show that the members of the
quantum $DS$-family are necessarily  compact quantum groups of Kac type and to obtain a description of their
`local' structure. Different parts of this `local' structure are combined in Section 6 to construct an example of
a compact quantum group which admits no non-trivial square root of the Haar state and yet is neither commutative
nor cocommutative. In that section we also state an open problem related to the possible `degree of complication'
of two-dimensional irreducible corepresentations of a quantum group in the quantum $DS$-family.

\section{Preliminaries}
\label{sec-prelim}

The symbol $\otimes$ will denote the spatial tensor product of $C^*$-algebras and $\odot$ the algebraic tensor
product, we use $\textup{Lin} F$ for the linear span of a set $F$ is a vector space and $\overline{\textup{Lin}}
F$
 for the closed linear span of a set $F$ in a Banach space.
\subsection{Compact quantum groups}

The notion of compact quantum groups has been introduced in
\cite{woronowicz87}. Here we adopt the definition from \cite{woronowicz98}
(Definition 2.1 of that paper).

\begin{definition}
A \emph{$C^*$-bialgebra} (a compact quantum semigroup) is a pair $(\mathsf{A}, \Delta)$, where
$\mathsf{A}$ is a unital $C^*$-algebra,
 $\Delta:\mathsf{A} \to \mathsf{A} \otimes \mathsf{A}$ is a unital,
 $*$-homomorphic map which is coassociative, i.e.\
\[ (\Delta \otimes \textup{id}_{\mathsf{A}})\circ \Delta = (\textup{id}_{\mathsf{A}} \otimes \Delta)
\circ\Delta.
\]
If the quantum cancellation properties
\[
\overline{{\rm Lin}}((1\otimes \mathsf{A})\Delta(\mathsf{A}) ) = \overline{{\rm Lin}}((\mathsf{A} \otimes
1)\Delta(\mathsf{A}) ) = \mathsf{A} \otimes \mathsf{A},
\]
are satisfied, then the pair $(\mathsf{A}, \Delta)$ is called a \emph{compact
  quantum group}.
\end{definition}

The map $\Delta$ is called the \emph{coproduct} of $\mathsf{A}$, it induces the convolution product
\[
\lambda\star \mu:=(\lambda\otimes \mu)\circ\Delta, \;\;\;
\lambda,\mu\in\mathsf{A}^*.
\]
When the coproduct is clear from the context we just speak of a compact quantum group $\mathsf{A}$.

The following fact is of the fundamental importance for this paper, cf.\ \cite[Theorem 2.3]{woronowicz98}.

\begin{proposition} \label{prop-haar}
Let $(\mathsf{A}, \Delta)$ be a compact quantum group. There exists a unique state $h \in \mathsf{A}^*$ (called the \emph{Haar state} of
$\mathsf{A}$) such that for all $a \in \mathsf{A}$
\[
(h \otimes \textup{id}_{\mathsf{A}})\circ  \Delta (a) = ( \textup{id}_{\mathsf{A}} \otimes h)\circ  \Delta (a) = h(a) 1.
\]
\end{proposition}

This naturally leads to the next definition introducing the main object of interest for the rest of the paper.

\begin{definition}
A state $\phi$ on a compact quantum group $\mathsf{A}$ is called a square root of the Haar state if
\[ \phi \star \phi = h.\]
It is said to be non-trivial if $\phi \neq h$.
\end{definition}

In general, the Haar state of a compact quantum group need not be faithful. But one can always divide by the nullspace of the
Haar state to produce a compact quantum group with faithful Haar state, usually called the reduced version of the original
quantum group, cf.\ \cite{bedos+murphy+tuset01}. This construction allows us to reduce our study to compact quantum groups with
faithful Haar states, see Lemma \ref{lem-faithful}.

\subsection{Quantum subgroups}

The notion of a quantum subgroup was introduced by Kac  \cite{kac68} in the setting of finite ring groups and by
Podle\'s \cite{podles95} for matrix pseudo-groups. In some contexts related to quantum subgroups it is necessary
to distinguish between the reduced and universal versions of the compact quantum groups in question (or consider
coamenable compact quantum groups, for which the two versions coincide), but it will not be important here.

\begin{definition}\label{def-subgroup}
A compact quantum group $(\mathsf{B},\Delta_{\mathsf{B}})$ is said to be a quantum subgroup of
a compact quantum group $(\mathsf{A},\Delta_{\mathsf{A}})$ if there exists a surjective compact quantum group morphism  $\pi:\mathsf{A} \to \mathsf{B}$, i.e.\ a surjective
unital
$*$-homomorphism $\pi:\mathsf{A} \to \mathsf{B}$ such that
\[
\Delta_{\mathsf{B}}\circ \pi = (\pi \otimes \pi) \circ\Delta_{\mathsf{A}}.
\]
A quantum subgroup $\mathsf{B}$ of $\mathsf{A}$ with Haar state $h_B$ is
called \emph{normal} if the images of the
conditional expectations
\begin{eqnarray*}
E_{\mathsf{A}/\mathsf{B}} &=& \big({\rm id}\otimes (h_{\mathsf{B}}\circ\pi)\big)\circ \Delta_{\mathsf{A}}, \\
E_{\mathsf{B}\backslash \mathsf{A}} &=& \big((h_{\mathsf{B}}\circ\pi)\otimes {\rm id}\big)\circ \Delta_{\mathsf{A}},
\end{eqnarray*}
coincide, cf.\ \cite[Proposition 2.1 and Definition 2.2]{wang08}. Note that the images of the conditional
expectations above can be thought of as the algebras of functions constant respectively on the right and left
`cosets' of the quantum subgroup $\mathsf{B}$.
\end{definition}

\subsection{Corepresentations}

An element $u =(u_{k\ell})_{1\le k,\ell\le n}\in M_n(\mathsf{A})$ is called an {\emph{$n$-dimensional
corepresentation of}} $(\mathsf{A}, \Delta)$ if for all $k,\ell=1,\ldots,n$ we have $\Delta(u_{k\ell}) =
\sum_{j=1}^n u_{kj} \otimes u_{j\ell}$.  All corepresentations considered in this paper are supposed to be
finite-dimensional. A corepresentation $u$ is said to be \emph{non-degenerate}, if $u$ is invertible,
\emph{unitary}, if $u$ is unitary, and \emph{irreducible}, if the only matrices $T\in M_n(\mathbb{C})$ with
$Tu=uT$ are multiples of the identity matrix. Two corepresentations $u,v\in M_n(\mathsf{A})$ are called
\emph{equivalent}, if there exists an invertible matrix $U\in M_n(\mathbb{C})$ such that $Uu=vU$.

An important feature of compact quantum groups is the existence of the dense $*$-subalgebra $\mathcal{A}$ (the
algebra of the \emph{smooth} elements of $\mathsf{A}$), which is in fact a Hopf $*$-algebra with the coproduct
$\Delta|_{\mathcal{A}}$ -- so for example $\Delta: \mathcal{A} \to \mathcal{A} \odot \mathcal{A}$. Fix a complete
family $(u^{(s)})_{s\in\mathcal{I}}$ of mutually inequivalent irreducible unitary corepresentations of
$(\mathsf{A}, \Delta)$, then $\{u^{(s)}_{k\ell}; s\in\mathcal{I},1\le k,\ell\le n_s\}$ (where $n_s$ denotes the
dimension of $u^{(s)}$) is a linear basis of $\mathcal{A}$, cf.\ \cite[Proposition 5.1]{woronowicz98}. We shall
reserve the index $s=\emptyset$ for the trivial representation $u^\emptyset=\mathbf{1}$.

Set $V_s={\rm span}\,\{u^{(s)}_{k\ell};1\le k,\ell\le n_s\}$ for $s\in\mathcal{A}$. By \cite[Proposition
5.2]{woronowicz98}, there exists a unique irreducible unitary corepresentation $u^{(s^{\rm c})}$, called the
\emph{contragredient} representation of $u^{(s)}$, such that $V_{s}^*=V_{s^{\rm c}}$. Clearly $(s^{\rm c})^{\rm
c}=s$.

The matrix elements of the irreducible unitary corepresentations satisfy the
famous Peter-Weyl orthogonality relations
\begin{equation}\label{eq-pw}
h\left(\left(u^{(s)}_{ij}\right)^* u^{(t)}_{k\ell}\right) = \frac{\delta_{st}
\delta_{j\ell} \overline{f\left(\left(u^{(s)}_{k
          i}\right)^*\right)}}{D_s}
\end{equation}
where $f:\mathcal{A}\to\mathbb{C}$ denotes the so-called Woronowicz character and
\[
D_s= \sum_{\ell=1}^{n_s} f\left(u_{\ell\ell}^{(s)}\right) = \sum_{\ell=1}^{n_s} \overline{f\left(\left(u_{\ell\ell}^{(s)}\right)^*\right)}
\]
is the quantum dimension of $u^{(s)}$, cf.\
\cite[Theorem 5.7.4]{woronowicz87}. Note that unitarity implies that the matrix
\[
f\left(\left(u^{(s)}_{k \ell}\right)^*\right)\in M_{n_s}(\mathbb{C})
\]
is invertible, with inverse $\big(f(u^{(s)}_{k\ell})\big)\in
M_{n_s}(\mathbb{C})$, cf.\ \cite[Equations (5.18), (5.24)]{woronowicz87}

We will say that a linear functional $\phi:\mathsf{A}\to\mathbb{C}$ has
\emph{finite support}, if
\[
\varphi|_{V_s}=0
\]
for all but finitely many $s\in\mathcal{I}$. The Haar state has
finite support, since it vanishes on all irreducible unitary representations
except the trivial one.

\begin{lemma}\label{lem-finite-support}
A continuous linear functional $\phi:\mathsf{A}\to\mathbb{C}$ has finite
support if and only if it admits a smooth density w.r.t.\ to the Haar state, i.e.\ if there
exists $x\in\mathcal{A}$ such that
\[
\phi(a) = h(xa), \qquad \mbox{ for all }a\in\mathsf{A}.
\]
The density $x$ is uniquely determined by $\phi$.
\end{lemma}
\begin{proof}
Assume such a density $x\in\mathcal{A}$ exists. We write $h_x$ for the linear
functional defined by $h_x(a)=h(xa)$ for all $a\in\mathsf{A}$. As a
smooth element, $x$ can be written as a finite linear combination
\[
x = \sum_{i=1}^n \sum_{k,\ell=1}^{n_{s_i}} c(s_i,k,\ell) u^{(s_i)}_{k\ell}
\]
Then the Peter-Weyl orthogonality relations \eqref{eq-pw} imply that
$\phi|_{V_s}=h_x|_{V_s}=0$ for $s\in \mathcal{I}$, $s\not\in\{s_1^{\rm
  c},\ldots,s_n^{\rm c}\}$.

Conversely, if $\phi$ has finite support, then the sum
\[
x= \sum_{s\in\mathcal{I}} \sum_{j,k,\ell=1}^{n_s} D_s \phi(u^{(s)}_{j\ell}) f(u^{(s)}_{kj}) \left(u^{(s)}_{k\ell}\right)^*
\]
is finite, therefore $x\in\mathcal{A}$, and the Peter-Weyl orthogonality
relations \eqref{eq-pw} imply
$\phi|_{\mathcal{A}}=h_x|_{\mathcal{A}}$. Density of $\mathcal{A}$ in
$\mathsf{A}$ and continuity of $\phi$ and $h_x$ then give $\phi=h_x$.

Clearly, by the Peter-Weyl orthogonality relations \eqref{eq-pw}, $x$ is
uniquely determined by $h_x|_{\mathcal{A}}=\phi|_{\mathcal{A}}$.
\end{proof}

\begin{remark}
Note that the density $x$ is uniquely determined by $\phi$, even if the Haar
state $h$ is not faithful. This is a consequence of the fact that the Haar
state is always faithful on the algebra $\mathcal{A}$ of smooth elements.
\end{remark}

Denote by $(\pi_h,H,\mathbf{1}_h)$ the GNS representation of $\mathsf{A}$ with respect to the Haar state. If the Haar state $h$
is faithful, we can make use of the Tomita-Takesaki theory for Haar states on compact quantum groups
\cite{woronowicz87,woronowicz98}. Define an antilinear operator $S_h$ on $H$ by
\[
S_h\pi_h(a)\mathbf{1}_h = \pi_h(a)^*\mathbf{1}_h,
\]
for any $a\in\mathsf{A}$ and set $\Delta_h=S^*_hS_h$. The modular automorphism group $(\sigma_t^h)_{t\in\mathbb{R}}$ is given by
\[
\pi_h\big(\sigma_t^h(a)\big)\mathbf{1}_h = \Delta_h^{it}\pi_h(a)\mathbf{1}_h
\]
for $a\in\mathcal{A}$. Each element of $\mathcal{A}$ is analytic with respect to the modular group
$(\sigma_t^h)_{t\in\mathbb{R}}$.

\subsection{Discrete quantum groups}\label{subsec-discrete}

Let $(\mathsf{A},\Delta)$ be a compact quantum groups. The space of linear functionals on $\mathsf{A}$ with
finite support has the structure of a discrete algebraic quantum group.

Fix a complete family $(u^{(s)})_{s\in\mathcal{I}}$ of mutually inequivalent irreducible
unitary corepresentations of $\mathsf{A}$, and define
$e^{(s)}_{k\ell}:\mathcal{A}\to\mathbb{C}$ for $s\in\mathcal{I}$, $1\le
k,\ell\le n_s$ by
\[
e^{(s)}_{k\ell}\left(u^{(t)}_{ij}\right) = \delta_{st} \delta_{ki}\delta_{\ell j}
\]
for $t\in\mathcal{I}$, $1\le i,j\le n_s$. These functionals extend to
continuous functionals on $\mathsf{A}$, since $e^{(s)}_{k\ell}=h_x$, with
$x=D_s  \sum_{j=1}^{n_s} f(u^{(s)}_{j\ell}) \left(u^{(s)}_{jk}\right)^*$.
The convolution product of two such functionals gives
$e^{(s)}_{ij}e^{(t)}_{k\ell}=\delta_{st} \delta_{jk} e^{(s)}_{i\ell}$
for $s,t\in\mathcal{I}$, $1\le i,j\le n_s$, $1\le k,\ell\le n_t$, i.e.\ the
linear functionals on $\mathsf{A}$ with finite support from a subalgebra
\[
\hat{\mathcal{A}}={\rm span}\left\{e^{(s)}_{ij};s\in\mathcal{I}, 1\le i,j\le n_s\right\}
\]
of $\mathsf{A}^*$ with respect to the convolution product. Equip $\hat{\mathcal{A}}$ with the involution
$\left(e^{(s)}_{k\ell}\right)^*=e^{(s)}_{\ell k}$. The $^*$-algebra $\hat{\mathcal{A}}$ has the form of a
multimatrix algebra,
\[
\hat{\mathcal{A}}=\bigoplus_{s\in\mathcal{I}}{\rm span}\left\{e^{(s)}_{\ell k};1\le k,\ell\le
    n_s\right\}\cong \bigoplus_{s\in\mathcal{I}}
    M_{n_s}(\mathbb{C})\subseteq\mathsf{A}^*
\]
(algebraic direct sum). With the coproduct $\hat{\Delta}:\hat{\mathcal{A}}\to
M(\hat{\mathcal{A}}\odot\hat{\mathcal{A}})$ defined by $\hat{\Delta}(\phi)(a\otimes b)=\phi(ab)$ for
$a,b\in\mathsf{A}$, $\hat{\mathcal{A}}$ becomes a discrete algebraic quantum group in the sense of
\cite{vandaele98,vandaele03}. Here $M(\hat{\mathcal{A}}\odot\hat{\mathcal{A}})$ denotes the multiplier algebra of
$\hat{\mathcal{A}}\odot\hat{\mathcal{A}}$, its elements can be naturally identified with linear functionals on
$\mathcal{A}\odot\mathcal{A}$.

For $J\subseteq \mathcal{I}$, we introduce the notation
\[
\hat{V}_J =
\left\{\phi\in\hat{\mathcal{A}};\phi\left(u^{(s)}_{ij}\right)=0\mbox{ for
  }s\not\in J,1\le i,j\le n_s\right\}
\]
for the space of functionals which vanish on the irreducible unitary corepresentations that do not belong to $J$.
We have
\[
\hat{V}_{\{s\}} = {\rm span}\left\{e^{(s)}_{\ell k};1\le k,\ell\le
    n_s\right\}\cong  M_{n_s}(\mathbb{C})
\]
for $s\in \mathcal{I}$ and
\[
\hat{V}_{\{s,s^c\}} \cong M_{n_s}(\mathbb{C})\oplus M_{n_s}(\mathbb{C})
\]
if $s\not=s^c$, i.e.\ if $u^{(s)}$ is not contragredient to itself.

The pair $(\hat{\mathcal{A}},\hat{\Delta})$ admits an antipode $\hat{S}:\hat{\mathcal{A}}\to\hat{\mathcal{A}}$,
which can be characterised by $\big(\hat{S}(\phi)\big)(a)=\phi\big(S(a)\big)$ for $a\in\mathcal{A}$,
$\phi\in\hat{\mathcal{A}}$.

The antilinear map $\mathsf{A}^*\ni\phi\mapsto \phi^\dag=\overline{\phantom{a}}\circ \phi\circ *\in\mathsf{A}^*$
allows to characterise the real algebra of hermitian functionals on $\mathsf{A}$ as its fixed point algebra. We
have $(\phi\star\psi)^\dag=\phi^\dag\star \psi^\dag$ and $(\phi^\dag)^\dag=\phi$. Since finitely supported
functionals are in the domain of the antipode $\hat{S}$, we have
$\phi^\dag=\hat{S}^{-1}(\phi^*)=\big(\hat{S}(\phi)\big)^*$ for $\phi\in\hat{\mathcal{A}}$.

\subsection{First examples}

\subsubsection{Commutative examples}

If $G$ is a compact group, then $\mathsf{A}=C(G)$ becomes a compact quantum
group with the coproduct $\Delta :\mathsf{A}=C(G)\to
\mathsf{A}\otimes\mathsf{A}\cong C(G\times G)$ defined by
\[
\Delta(f)(g_1,g_2)= f(g_1g_2),
\]
for $f\in C(G)$, $g_1,g_2\in G$. Furthermore, any commutative compact quantum
group is of this form, cf.\ \cite[Remark 3 following Definition
1.1]{woronowicz98}.

The Haar state on a commutative compact quantum group $C(G)$ is given by
integration against the Haar measure $\mu$ of $G$, i.e.\ $h(f)=\int_G f{\rm
  d}\mu$ for $f\in C(G)$. It admits a non-trivial square root if and only if
the Haar measure $\mu$ admits a non-trivial square root. Hence the main theorem of \cite{diaconis+shahshahani86}
can be reformulated in the following way.

\begin{theorem}[\cite{diaconis+shahshahani86}] \label{DS}
Let $G$ be a separable compact group. The pair $(C(G), \Delta))$ admits no non-trivial square root of the Haar
state if and only if $G$ is abelian or of the form $H \times E$ where $H$ is the group of unit quaternions and
$E$ is a Cartesian product of (at most countably many) copies of $\mathbb{Z}_2$.
\end{theorem}

\subsubsection{Cocommutative examples}

A compact quantum group $(\mathsf{A},\Delta)$ is called \emph{cocommutative},
if $\tau\circ\Delta=\Delta$, where
$\tau:\mathsf{A}\otimes\mathsf{A}\to\mathsf{A}\otimes\mathsf{A}$ is the flip,
$\tau(a\otimes b)=b\otimes a$. Since irreducible corepresentations of a
cocommutative compact quantum group are necessarily one-dimensional, the dense $*$-Hopf algebra of smooth elements in
a cocommutative compact quantum group $(\mathsf{A},\Delta)$ is of the form
$\mathcal{A}={\rm span}(\Gamma)$, where
\[
\Gamma = \{u\in \mathsf{A};u\mbox{ unitary and }\Delta(u)=u\otimes u\}
\]
is a (discrete) subgroup of the group $\mathcal{U}(\mathsf{A})$ of unitary
elements of $\mathsf{A}$.

The Haar
state $h$ acts as
\[
h(u)=\left\{\begin{array}{cll}
1 & \mbox{ if } u=\mathbf{1} \mbox{ the trivial corepresentation}, \\
0 & \mbox{ else.}
\end{array}\right.
\]
on $u\in\Gamma$.
If $\phi$ is a square root of $h$, then
$(\phi\star\phi)(u)=\big(\phi(u)\big)^2=\delta_{u\mathbf{1}}$ for
$u\in\Gamma$, and the only positive square root of the Haar state is the
trivial solution $\phi=h$. This argument shows that the Haar state of a cocommutative
compact quantum group never admits a non-trivial square root.

\subsection{Terminology}

Motivated by Theorem \ref{DS} we introduce the following terminology.

\begin{definition} \label{DSfamily}
A compact quantum group $\mathsf{A}$ is said to belong to the quantum $DS$-family if its Haar state does not
admit any non-trivial square roots.
\end{definition}

We will sometimes refer to a \emph{classical $DS$-family} as the family of groups listed in Theorem \ref{DS}. The examples discussed above show that one can find and characterise commutative and cocommutative quantum groups in the quantum $DS$-family. A priori it is not clear if there exist at all any elements in the quantum $DS$-family which belong to neither of these classes. We will in fact exhibit such an example in Proposition \ref{prop-quantexam}.

\section{Square roots of the Haar state and nilpotent functionals}
\label{sec-main}

We will show in this section that a compact quantum group $(\mathsf{A},\Delta)$ is not in the quantum $DS$-family
if and only if there exists a hermitian functional on $\mathsf{A}$ that is nilpotent for the convolution, cf.\
Theorem \ref{thm-nilp}. Necessity of this condition is immediate. Indeed, if a state
$\phi:\mathsf{A}\to\mathbb{C}$, $\phi\not=h$, is a square root of the Haar state, then $\rho=\phi-h\not=0$ is
hermitian and nilpotent, since
\[
\rho\star\rho=\phi\star \phi - \phi\star h - h\star \phi + h\star h =0.
\]
To prove the converse, we follow a similar strategy as Diaconis and
Shahshahani. Given a nilpotent hermitian functional $\rho$, we construct a new
(``truncated'') hermitian functional $\psi_s$ for which there exists
$\epsilon\not=0$ such that $h+\epsilon \psi_s$ defines a state which is a square
root of the Haar state $h$. But if the Haar state is not a trace, then more
care is required to prove the positivity of $h+\epsilon \psi_s$.

We begin with a simple lemma for the tracial case.
\begin{lemma}\label{lem-boundedness-tracial}
Let $\mathsf{A}$ be a unital $C^*$-algebra with tracial state $h$. Then we have
\[
|h(xa)| \le ||x|| h(a)
\]
for all $x\in \mathsf{A}$ and $a\in \mathsf{A}_+$.
\end{lemma}

\begin{proof}
Since $a$ is positive, there exists $b\in\mathsf{A}$ such that $a=bb^*$. Denote by $\pi_h$ and $\mathbf{1}_h$ the
GNS representation of $h$ and its cyclic vector representing the state $h$. Then we have
\begin{eqnarray*}
|h(xa)| &=& |h(xbb^*)| = |h(b^*xb)| = \langle \pi_h(b)
\mathbf{1}_h,\pi_h(x)\pi_h(b)\mathbf{1}_h\rangle \\
&\le & ||x|| \langle \pi_h(b)
\mathbf{1}_h,\pi_h(b)\mathbf{1}_h\rangle = ||x|| h(b^*b) = ||x|| h(a),
\end{eqnarray*}
since $\pi_h$ is a contraction.
\end{proof}

Let us now characterise hermitianity and positivity of a given finitely supported linear functional in terms of its density.
Recall that elements in $\mathcal{A}$ are analytic with respect to the modular automorphism group of the Haar state.

\begin{lemma}\label{lem-hermitian-pos}
Let $(\mathsf{A},\Delta)$ be a compact quantum group with faithful Haar state $h$ and modular group
$(\sigma^h_t)_{t \in \mathbb{R}}$, and let $x\in\mathcal{A}$.
\begin{enumerate}
\item
The functional $h_x\in \mathsf{A}^*$, $h_x(a)=h(xa)$ is hermitian if and only
if $\sigma^h_{-i/2}(x)$ is self-adjoint.
\item
The functional $h_x\in \mathsf{A}^*$, $h_x(a)=h(xa)$ is positive if and only
if $\sigma^h_{-i/2}(x)$ is positive.
\end{enumerate}
\end{lemma}

\begin{proof}
\begin{enumerate}
\item Denote by $\mathsf{A}^*_h$ the space of hermitian continuous functionals on $\mathsf{A}$ and once again
write  $\pi_h$ and $\mathbf{1}_h$ for the GNS representation of $h$ and its cyclic vector representing the state
$h$. We have
\begin{eqnarray*}
h_x\in \mathsf{A}^*_h & \Leftrightarrow & h(xa^*)= \overline{h(xa)} \quad\forall
a\in\mathsf{A} \quad(\mbox{or } a\in\mathcal{A})\\
&\Leftrightarrow& \langle x^*\mathbf{1}_h,a^*\mathbf{1}_h\rangle =\overline{\langle x^*\mathbf{1}_h,a\mathbf{1}_h\rangle} \quad\forall
a\in\mathcal{A} \\
&\Leftrightarrow& \langle S_h x\mathbf{1}_h,S_ha\mathbf{1}_h\rangle =\langle a\mathbf{1}_h,x^*\mathbf{1}_h\rangle \quad\forall
a\in\mathcal{A} \\
&\Leftrightarrow& \langle a\mathbf{1}_h,\underbrace{\Delta_h x\mathbf{1}_h}\rangle =\langle a\mathbf{1}_h,x^*\mathbf{1}_h\rangle \quad\forall
a\in\mathcal{A} \\
&& \hspace{16mm}= \sigma^h_{-i}(x)\mathbf{1}_h \\
&\Leftrightarrow& \sigma_{-i}^h(x)=x^* \quad\Leftrightarrow\quad \sigma^h_{-i/2}(x) =
\big(\sigma^h_{-i/2}(x)\big)^*,
\end{eqnarray*}
where we used faithfulness of $h$ and the relation $\sigma^h_{i/2}\circ *=*\circ\sigma^h_{-i/2}$.
\item
Denote by $\mathsf{A}^*_+$ the space of positive functionals on
$\mathsf{A}$. We have
\begin{eqnarray*}
h_x\in\mathsf{A}^*_+ & \Leftrightarrow & h(xbb^*) \ge 0 \quad \forall
b\in\mathsf{A} \quad(\mbox{or }b\in \mathcal{A}) \\
& \Leftrightarrow & h\big(\sigma_i^h(b^*)xb\big) \ge 0 \quad \forall
b\in \mathcal{A} \\
& \Leftrightarrow & \langle \sigma^h_{-i}(b)\mathbf{1}_h,xb\mathbf{1}_h\rangle \ge 0 \quad \forall b\in
\mathcal{A}.
\end{eqnarray*}
Since
\begin{eqnarray*}
\langle \sigma^h_{-i}(b)\mathbf{1}_h,xb\mathbf{1}_h\rangle &=& \langle
\Delta_hb\mathbf{1}_h,xb\mathbf{1}_h\rangle = \langle
\Delta_h^{1/2}b\mathbf{1}_h,\Delta_h^{1/2}xb\mathbf{1}_h\rangle \\
&=& \langle
\Delta_h^{1/2}b\mathbf{1}_h,\sigma_{-i/2}^h(x)\Delta_h^{1/2}b\mathbf{1}_h\rangle
\end{eqnarray*}
and since
$\{\Delta_h^{1/2}b\mathbf{1}_h;b\in\mathcal{A}\}=\{\sigma_{-i/2}^h(b)\mathbf{1}_h;b\in\mathcal{A}\}$
is dense, this is equivalent to $\sigma_{-i/2}^h(x)\ge 0$.
\end{enumerate}
\end{proof}

\begin{lemma}\label{lem-pos}
Let $x\in\mathcal{A}$. If the functional $h_x\in \mathsf{A}^*$ is hermitian, then there exists $\epsilon>0$ such that $h+\epsilon
h_x$ is a positive.
\end{lemma}

\begin{proof}
Set $\varphi_\epsilon=h+\epsilon h_x=h_{\mathbf{1}+\epsilon x}$. Since $h_x\in
\mathsf{A}^*_h$, $\sigma_{-i/2}^h(x)$ is self-adjoint by Lemma
\ref{lem-hermitian-pos}. Therefore there exists $\epsilon>0$ such that
$\mathbf{1}\ge \epsilon \sigma_{-i/2}^h(x)\ge -\mathbf{1}$. Then
$\sigma_{-i/2}^h(\mathbf{1}+\epsilon x)= \mathbf{1}+\epsilon \sigma_{-i/2}^h(x)\ge
\mathbf{1}-\mathbf{1}=0$. Since $\mathbf{1}+\epsilon x\in\mathcal{A}$, we can
apply Lemma \ref{lem-hermitian-pos} and get $\varphi_\epsilon\in\mathsf{A}^*_+$.
\end{proof}

\begin{remark}
Similar methods yield the following general result.

If $M$ is a von Neumann algebra with a faithful normal state $\omega$ and $x \in M$ analytic with respect to the modular
automorphism group $\{\sigma_t: t \in \mathbb{R}\}$ of the state $\omega$, then
\[
|\omega(xa)| \leq \|\sigma_{-\frac{i}{2}}(x)\| \omega(a), \qquad \mbox{ for }a \in M_+.
\]
\end{remark}

\begin{lemma}\label{lem-trunc}
\emph{(Truncation Lemma)}
Let $\rho\in\mathsf{A}^*$ be a hermitian functional such that
$\rho\star\rho=0$. For $u^{(s)}$ an irreducible unitary corepresentation of
$\mathsf{A}$, define $\psi_s$ by
\[
\psi_s\left(u^{(t)}_{k\ell}\right) = \left\{
\begin{array}{cl}
\rho\left(u^{(t)}_{k\ell}\right) & \mbox{ if } t\mbox{ is equivalent to
}s\mbox{ or }s^c, \\
0 & \mbox{ else.}
\end{array}\right.
\]
Then $\psi_s$ is hermitian, has finite support, and satisfies $\psi_s\star\psi_s=0$.
\end{lemma}

\begin{proof}
The support of $\psi_s$ is contained in the $*$-closed subspace $V_s+V_{s^c}$
and $\psi_s|_{V_s+V_{s^c}}=\rho|_{V_s+V_{s^c}}$. Therefore $\psi_s$ is clearly
hermitian and finitely supported. Furthermore,
\begin{gather*}
(\psi_{s}\star\psi_s)(u^{(t)}_{k\ell})=\sum_{j=1}^{n_t}
\psi_s(u_{kj}^{(t)})\psi_s(u_{j\ell}^{(t)}) \\
= \left\{\begin{array}{ll} \sum_{j=1}^{n_t}
\rho(u_{kj}^{(t)})\rho(u_{j\ell}^{(t)})
=(\rho\star\rho)(u^{(t)}_{k\ell})=0\qquad&\mbox{ if }t\mbox{ is equivalent to }s\mbox{ or }s^c, \\
0 \qquad &\mbox{ else}
\end{array}\right.
\end{gather*}
for all irreducible unitary corepresentations $u^{(t)}$ of $\mathsf{A}$, i.e.\ $\psi_s\star\psi_s=0$.
\end{proof}

Let us first show that it is sufficient to consider compact quantum groups
with faithful Haar states.

\begin{lemma}\label{lem-faithful}
Let $(\mathsf{A},\Delta)$ be a compact quantum groups with not necessarily
faithful Haar state $h$ and denote by $(\tilde{\mathsf{A}},\tilde{\Delta})$ its reduced
version, i.e.\ the compact
quantum group with faithful Haar state $\tilde{h}$ obtained from
$(\mathsf{A},\Delta)$ by dividing out the nullspace of $h$.

If $\tilde{h}$ admits a non-trivial square root, then so does $h$.
\end{lemma}

\begin{remark}
The converse is also true, and can be shown using truncation arguments similar to those in the proof of Theorem \ref{thm-nilp},
but we will not need it.
\end{remark}

\begin{proof}
Denote by $\tilde{\pi}:\mathsf{A}\to\tilde{\mathsf{A}}$ the canonical
projection from $(\mathsf{A},\Delta)$ to $(\tilde{\mathsf{A}},\tilde{\Delta})$, cf.\ \cite{bedos+murphy+tuset01}. Then
$h=\tilde{h}\circ\tilde{\pi}$ and if $\tilde{h}$ admits a non-trivial square
root $\tilde{\phi}\not=h$, then clearly
$\phi=\tilde{\phi}\circ\tilde{\pi}\not=h$ defines a non-trivial square root of
$h$.
\end{proof}

\begin{theorem}\label{thm-nilp}
The Haar state $h$ of a compact quantum group $(\mathsf{A},\Delta)$ admits a non-trivial square root, i.e.\ a
state $\phi\not=h$ such that $\phi\star\phi=h$, if and only if there exists a bounded non-zero hermitian
continuous linear functional on $\mathsf{A}$ that is nilpotent for the convolution product.
\end{theorem}

\begin{proof}
If $h$ admits a non-trivial square root $\phi$, then clearly $\rho=\phi-h$ defines a bounded non-zero hermitian
nilpotent functional.

Conversely, assume that $\mathsf{A}^*$ contains a non-zero nilpotent hermitian
functional $\psi$. Then all convolution powers of $\psi$ are
also hermitian. If $\psi^{\star n}=0$, and $n$ is the smallest such number, then
set $\rho=\psi^{\star (n-1)}$. This is non-zero and satisfies
$\rho\star\rho=0$. Therefore $\rho(u)=0$ for any $u\in\mathsf{A}$ with
$\Delta(u)=u\otimes u$, in particular $\rho(\mathbf{1})=0$.

Since $\rho\not=0$, there exists an irreducible unitary corepresentation
$u^{(s)}$ such that $\rho|_{V_s+V_{s^c}}\not=0$. Fix such an irreducible
unitary corepresentation $u^{(s)}$ and define $\psi_s$ as in the Truncation Lemma
(Lemma \ref{lem-trunc}). Then $\psi_s$ has finite support and there exists a
unique $x\in\mathcal{A}$ such that $\psi_s=h_x\in\mathsf{A}^*_h$, cf.\ Lemma \ref{lem-finite-support}.

If $h$ is faithful, then, by Lemma
\ref{lem-pos}, there exists $\epsilon>0$ such that $\phi=h+\epsilon h_x\in
\mathsf{A}_+^*$. Since $h\star \psi_s=\psi_s\star h = \psi_s(\mathbf{1})h =0$,
we get
\[
\phi\star\phi = h\star h + \epsilon h\star\psi_s + \epsilon \psi_s\star h +
\epsilon^2 \psi_s\star\psi_s = h,
\]
i.e.\ $\phi$ is a non-trivial square root of the Haar state $h$.

If $h$ is not faithful, then $x\in\mathcal{A}$ can be used to define a nilpotent hermitian
functional $\tilde{\psi_s}=\tilde{h}_x$ and a non-trivial square
$\tilde{\phi}=\tilde{h}+\epsilon \tilde{h}_x$ on the reduced version
$(\tilde{\mathsf{A}},\tilde{\Delta})$.
By Lemma \ref{lem-faithful}, $\phi=\tilde{\phi}\circ\tilde{\pi}$ then defines a
non-trivial square root of $h$ on $(\mathsf{A},\Delta)$.
\end{proof}

Since by the Truncation Lemma we can always choose this nilpotent hermitian
linear functional to have finite support, we also get the following
characterisation.

\begin{corollary}\label{cor-nilp}
The Haar state $h$ of a compact quantum group $(\mathsf{A},\Delta)$ admits a non-trivial square root if and only
if its dual discrete algebraic quantum group $(\hat{\mathcal{A}},\hat{\Delta})$ introduced in Subsection
\ref{subsec-discrete} contains a non-zero nilpotent element that is hermitian w.r.t.\ $\dag$.
\end{corollary}

We will now show that compact quantum group whose Haar state has no non-trivial square root is hamiltonian and
that all its quantum subgroups are normal.

\begin{lemma}\label{lem-normal}
Let $(\mathsf{A},\Delta)$ be a compact quantum group with faithful Haar state. A quantum subgroup $\mathsf{B}$ of
$\mathsf{A}$ is normal if and only if the idempotent state $h_{\mathsf{B}}\circ\pi$ on $\mathsf{A}$ induced by
the Haar state $h_{\mathsf{B}}$ of $\mathsf{B}$ is in the center of $\mathsf{A}^*$.
\end{lemma}

\begin{proof}
Denote by
\begin{eqnarray*}
E_{\mathsf{A}/\mathsf{B}} &=& \big({\rm id}\otimes (h_{\mathsf{B}}\circ\pi)\big)\circ \Delta, \\
E_{\mathsf{B}\backslash \mathsf{A}} &=& \big((h_{\mathsf{B}}\circ\pi)\otimes {\rm id}\big)\circ \Delta,
\end{eqnarray*}
the conditional expectations onto the coidalgebras $\mathsf{A}/\mathsf{B}$ and $\mathsf{B}\backslash\mathsf{A}$.
The quantum subgroup $\mathsf{B}$ is normal if and only these two coidalgebras coinside, cf.\ Definition
\ref{def-subgroup}. Since $E_{\mathsf{A}/\mathsf{B}}$ and $E_{\mathsf{B}\backslash \mathsf{A}}$ are unital and
preserve the Haar state, by the uniqueness of state preserving conditional expectations this
  is equivalent to $E_{\mathsf{A}/\mathsf{B}}=E_{\mathsf{B}\backslash \mathsf{A}}$, or
\[
f\star (h_{\mathsf{B}}\circ\pi) = f\circ E_{\mathsf{A}/\mathsf{B}} = f\circ
E_{\mathsf{B}\backslash \mathsf{A}}= (h_{\mathsf{B}}\circ \pi)\star f
\]
for all $f\in \mathsf{A}^*$.
\end{proof}

\begin{definition}
We call a compact quantum group $(\mathsf{A},\Delta)$ \emph{hamiltonian} if all idempotent states on $\mathsf{A}$
are central in $\mathsf{A}^*$ (w.r.t.\ the convolution).
\end{definition}

Idempotent states on finite and compact quantum groups were characterised in
\cite{franz+skalski-idemp,franz+skalski-note09}. For a compact group $G$, all idempotent states on $C(G)$ are
induced by Haar measures of closed subgroups of $G$, and $C(G)$ is hamiltonian if and only if all closed
subgroups of $G$ are normal. Lemma \ref{lem-normal} shows that all quantum subgroups of hamiltonian compact
quantum groups have to be normal. But noncommutative compact quantum groups may have idempotent states that are
not induced from quantum subgroups.

\begin{proposition}\label{prop-normal}
Let $(\mathsf{A},\Delta)$ be a compact quantum group in the quantum $DS$-family. If the Haar state of $\mathsf{A}$
is faithful,  then $(\mathsf{A}, \Delta)$ is hamiltonian. In particular, every quantum subgroup of $(\mathsf{A}, \Delta)$ is normal.
\end{proposition}
\begin{proof}
By Theorem \ref{thm-nilp}, if $h_A$ admits no non-trivial square root,
then $\mathsf{A}^*$ contains no hermitian functionals that are nilpotent for
the convolution product.

Then the result follows from the fact that
in a unital ring without nilpotent elements all idempotents are central, as in \cite[Lemma 3]{diaconis+shahshahani86}. Since $\mathsf{A}_h^*=\{\phi\in \mathsf{A}^*; \phi\mbox{ hermitian}\}$ has no
nilpotent elements, all idempotent states are central in $\mathsf{A}_h^*$, and therefore also in $\mathsf{A}^*$.
\end{proof}

As mentioned in the introduction the classification of  hamiltonian groups plays a very important role in the arguments of
\cite{diaconis+shahshahani86}. As no such classification is known for (compact) quantum groups, to study compact quantum groups
which do not admit nontrivial square roots of Haar states we need to develop other techniques. The next two sections will be
devoted to this task.

\section{A structure theorem}
\label{sec-structure}

In this section we characterise compact quantum groups whose Haar state admits no non-trivial square root in
terms of their irreducible unitary corepresentations or their dual discrete algebraic quantum group, see Theorem
\ref{thm-structure} and Proposition \ref{prop-char}.

By Corollary \ref{cor-nilp}, we have to characterise discrete algebraic
quantum groups which contain no non-zero hermitian nilpotent elements.

For $s\in\mathcal{I}$, we define $R_s$ to be the real algebra of
hermitian linear functionals on $\mathsf{A}$ that vanish on all irreducible
unitary corepresentations of $(\mathsf{A},\Delta)$ except $s$ and $s^c$, i.e.
\[
R_s= \hat{V}_{\{s,s^c\}}\cap \mathsf{A}^*_h.
\]
Clearly $R_s= R_{s^c}$.
Since $\hat{\mathcal{A}}$ is a multimatrix algebra, the real algebra of all finitely
supported hermitian linear functionals $\hat{\mathcal{A}}_h =
\hat{\mathcal{A}}\cap \mathsf{A}^*_h=$ decomposes into a direct sum
\begin{equation}\label{eq-decomp}
\hat{\mathcal{A}}_h = \bigoplus_{s\in\mathcal{I}_r} R_s,
\end{equation}
where the direct sum runs over the reduced index set $\mathcal{I}_r$ which is
obtained from $\mathcal{I}$ by choosing only one representative from each set $\{s,s^c\}$.

Recall that Frobenius \cite{frobenius1878} has shown that there exist exactly
three finite-dimensional division algebras over $\mathbb{R}$, namely the field of real numbers
$\mathbb{R}$, the field of complex numbers $\mathbb{C}$, and the skew field of
quaternions $\mathbb{H}$.

\begin{theorem}\label{thm-structure}
A compact quantum group $(\mathsf{A},\Delta)$ belongs to the quantum $DS$-family if and only if all summands occuring in the
decomposition \eqref{eq-decomp} are isomorphic to one of the three finite-dimensional division algebras $\mathbb{R}$, $\mathbb{C}$, or $\mathbb{H}$.
\end{theorem}
\begin{proof}
Let us first verify that the conditions are sufficient. If all summands $R_s$ in the direct sum \eqref{eq-decomp} are isomorphic
to division algebras, then they can not contain non-zero nilpotent elements, therefore $\hat{\mathcal{A}}_h$ has no non-zero
nilpotent elements, either, and by Corollary \ref{cor-nilp} the Haar state of $(\mathsf{A},\Delta)$ admits no non-trivial square
root.

Conversely, if the Haar state of $(\mathsf{A},\Delta)$ has no non-trivial square root, then none of the real
algebras $R_s$ occuring in \eqref{eq-decomp} contain non-zero nilpotent elements.

Let $s\in\mathcal{I}$. If $u^{(s)}$ is not contragredient to itself, i.e.\ if
$u^{(s)}$ and $u^{(s^c)}$ are not equivalent, then we have
$\hat{V}_{\{s\}}\cong M_{n_s}(\mathbb{C})$, $V_s\cap V_{s^c}=\{0\}$, and the map
\[
\hat{V}_{\{s\}}\ni\phi\mapsto \phi+\phi^\dag\in R_s\subseteq \hat{V}_{\{s,s^c\}},
\]
where
\[
(\phi+\phi^\dag)(a) = \left\{\begin{array}{ll}
\phi(a) &\mbox{ if }a\in V_s, \\[3pt]
\phi^\dag(a)=\overline{\phi(a^*)} &\mbox{ if }a\in V_{s^c}, \\[1pt]
0 &\mbox{ if }a\in V_t\mbox{ with }t\not=s,s^c,
\end{array}\right.
\]
for $\phi\in\hat{V}_{\{s\}}$ is an isomorphism of real algebras. Therefore
$R_s\cong\hat{V}_{\{s\}}\cong M_{n_s}(\mathbb{C})$ as real algebras. We see
that $R_s$ contains no non-zero nilpotent elements if and only if $n_s=1$, and in
this case $R_s\cong \mathbb{C}$.

Let us now consider the case where $u^{(s)}$ is contragredient to itself. Then
we have $V_s=V_{s^c}$, $\hat{V}_{\{s\}}=\hat{V}_{\{s,s^c\}}\cong
M_{n_s}(\mathbb{C})$. $R_s$ is a real subalgebra of
$\hat{V}_{\{s\}}$ whose complexification coincides
with $\hat{V}_{\{s\}}$, since any linear functional on $V_s$ can be written as
a complex linear combination of hermitian functionals. $\hat{V}_{\{s\}}\cong
M_{n_s}(\mathbb{C})$ is simple, and since the complexification of any real
ideal in $R_s$ would be an ideal in $\hat{V}_{\{s\}}$, it follows that $R_s$
is also simple. Therefore, by Wedderburn's theorem (\cite{wedderburn08}, see
also \cite[Section 13.11]{vanderwaerden91}), $R_s$ is isomorphic to a full matrix algebra over a division
algebra. In other words, we have $R_s\cong M_m(\mathbb{K})$ for some $m\ge 1$
and $\mathbb{K}\in\{\mathbb{R},\mathbb{C},\mathbb{H}\}$. If $R_s$ contains no
non-zero nilpotent elements, then necessarily $m=1$, and only the two cases
$n_s=1$ and $R_s\cong \mathbb{R}$, or $n_s=2$ and $R_s\cong\mathbb{H}$ can occur.
\end{proof}

Let us describe which corepresentations of a compact quantum group $(\mathsf{A},\Delta)$ lead to a real algebra of hermitian
functionals that is isomorphic to one of the three finite-dimensional division algebras $\mathbb{R}$, $\mathbb{C}$, and
$\mathbb{H}$.

Let $u=(u_{jk})_{1\le j,k\le n}\in M_n(\mathsf{A})$ be an irreducible unitary corepresentation of a
compact quantum group $(\mathsf{A},\Delta)$. We denote by $\bar{u}\in
M_n(\mathsf{A})$ the corepresentation obtained by taking the adjoints of the
coefficients of $u$, i.e.\ $\bar{u}= (u_{jk}^*)_{1\le j,k\le n}$. The
corepresentation $\bar{u}$ is not necessarily unitary, but it is
non-degenerate and equivalent to the contragredient corepresentation of $u$,
i.e.\ there exists an invertible matrix $T\in M_n(\mathbb{C})$ s.t.\
$\bar{u}=Tu^cT^{-1}$, cf.\ \cite[Proposition 5.2]{woronowicz98} or
\cite[Proposition 6.10]{maes+vandaele98}.

\begin{proposition}\label{prop-char}
Let $u\in M_n(\mathsf{A})$ be an irreducible unitary corepresentation of a compact quantum group
$(\mathsf{A},\Delta)$ and denote by $R(u)\subseteq \mathsf{A}^*_h$ the real algebra given by hermitian linear
functionals on $\mathsf{A}$ which vanish on the coefficients of all irreducible unitary corepresentations that
are not equivalent to $u$ or $u^c$. Then we have the following characterisations of $R(u)$.
\begin{itemize}
\item[(i)]
$R(u)\cong \mathbb{R}$ if and only if $u$ is one-dimensional and
contragredient to itself. This is the case if and only if $u$ is unitary,
self-adjoint, and group-like (i.e.\ $\Delta(u)=u\otimes u$).
\item[(ii)]
$R(u)\cong \mathbb{C}$ if and only if $u$ is one-dimensional and
not contragredient to itself. This is the case if and only if $u$ is unitary
and group-like, but not self-adjoint.
\item[(iii)]
$R(u)\cong \mathbb{H}$ if and only if $u$ is two-dimensional and there exists
an invertible matrix $Q\in M_2(\mathbb{C})$ such that $\bar{u}=TuT^{-1}$,
where $T=\bar{Q}\left(\begin{array}{cc} 0 & 1 \\ -1 & 0\end{array}\right)
Q^{-1}$. In particular, this implies that $u$ is contragredient to itself.
\end{itemize}
\end{proposition}
\begin{proof}
The first two cases follow from the proof of Theorem
\ref{thm-structure} if we note that a one-dimensional unitary corepresentation
is contragredient to itself if and only if it is self-adjoint.

Let us now  prove (iii). The real division algebra of
quaternions can be realised as
\[
\mathbb{H}=\left\{\left(\begin{array}{cc} \alpha + i\beta & i\gamma-\delta \\
      i\gamma+\delta & \alpha-i\beta\end{array}\right) : \alpha,\beta,\gamma,\delta\in\mathbb{R}\right\}.
\]
Its complexification is isomorphic to $M_2(\mathbb{C})$, and the elements of $\mathbb{H}$ can be characterised in
$M_2(\mathbb{C})$ as the hermitian elements for the anti-linear homomorphism $\dag:M_2(\mathbb{C})\to
M_2(\mathbb{C})$,
\[
\left( \begin{array}{cc} a & b \\ c & d\end{array}\right)^\dag = \left( \begin{array}{cc} \overline{d} & -\overline{c} \\
-\overline{b} & \overline{a}\end{array}\right).
\]
Dualising these relations we see that the real algebra $R(u)$ associated to a unitary
corepresentation $u$ is isomorphic to $\mathbb{H}$ if and only if $u$ is
two-dimensional and if the subspace $V(u)$ spanned by the coefficients of $u$
admits a basis $a_{11}, a_{12}, a_{21}, a_{22}$ such that
\[
a_{11}^* = a_{22}, \quad a_{12}^* =-a_{21}, \quad \mbox{ and
}\quad\Delta(a_{jk})=\sum_{\ell = 1}^2 a_{j\ell}\otimes a_{\ell k} \mbox{ for }1\le
  j,k\le 2.
\]
Therefore $a=\left(\begin{array}{cc} a_{11} & a_{12} \\ a_{21} &
    a_{22}\end{array}\right)\in M_2(\mathbb{C})$ is a corepresentation of
$(\mathsf{A},\Delta)$. It is non-degenerate, since its coefficients form a basis of $V(u)$, so by \cite[Proposition
6.4]{maes+vandaele98} it is equivalent to a unitary corepresentation, which we can choose to be $u$. I.e.\ there exists an
invertible matrix $Q\in M_2(\mathbb{C})$ s.t.\ $u=QaQ^{-1}$. We get
\[
\bar{u} = \bar{Q}\bar{a}\bar{Q}^{-1} = \bar{Q}Fa F^{-1}\bar{Q}^{-1}=TuT^{-1}
\]
with $F=\left(\begin{array}{cc} 0 & 1 \\ -1 & 0 \end{array}\right)$ and $T=\bar{Q}FQ^{-1}$.
\end{proof}

Note that the characterisation above seems to be new even for standard compact groups, together with Theorem \ref{DS} yielding the following corollary.

\begin{corollary} \label{corclassic}
Let $G$ be a separable compact group. The following conditions are equivalent:
\begin{rlist}
\item $G$ admits only one-dimensional and two-dimensional irreducible representations, each  two-dimensional irreducible representation $U:G \to M_2(\mathbb{C})$ is self-contragredient and there exists
an invertible matrix $Q\in M_2(\mathbb{C})$ such that $\bar{U}=TUT^{-1}$,
where $T=\bar{Q}\left(\begin{array}{cc} 0 & 1 \\ -1 & 0\end{array}\right)
Q^{-1}$;
\item $G \approx H \times E$ where $H$ is the group of unit quaternions and
$E$ is a Cartesian product of (at most countably many) copies of $\mathbb{Z}_2$.
\end{rlist}
\end{corollary}

Let us now consider the Woronowicz quantum group $SU_q(2)$. This example will play an important role in the next
section, when we show that a compact quantum group whose Haar state admits no non-trivial square root is
necessarily of Kac type, i.e.\ its Haar state is a trace, cf.\ Theorem \ref{thm-kac-type}.

\begin{example}\label{exa-suq}
Let $q\in \mathbb{R}\backslash\{0\}$.
Denote by $SU_q(2)=(\mathsf{A},\Delta)$ the Woronowicz quantum group \cite{woronowicz87b}, i.e.\
the universal $C^*$-algebra generated by the four generators
$u_{11},u_{12},u_{21},u_{22}$ with the coproduct determined by $\Delta u_{jk}
= \sum_{\ell=1}^2 u_{j\ell}\otimes u_{\ell k}$ for $j,k=1,2$, and the $*$-algebraic relations $uu^*=I=u^*u$
and $\bar{u}=F_quF_q^{-1}$ with $u=\left(\begin{array}{cc} u_{11} & u_{12} \\
    u_{21} & u_{22}\end{array}\right)$ and $F_q=\left(\begin{array}{cc} 0 & q \\
    -1 & 0\end{array}\right)$, i.e.\
\begin{eqnarray*}
\sum_{\ell =1}^2 u_{j\ell} u^*_{k \ell} = &\delta_{jk}& = \sum_{\ell =1}^2
u^*_{\ell j} u_{\ell k}, \\
u_{11}^* = u_{22}, &\quad& u_{12}^* = -q u_{21}.
\end{eqnarray*}
Note that $SU_q(2)$ is isomorphic to the universal orthogonal quantum group
$A_o(\tilde{F}_q)$ defined by Wang and van Daele \cite{vandaele+wang96}, where
$\tilde{F}_q$ is given by
\[
\tilde{F}_q = \frac{F_q}{\sqrt{|\det F_q|}}=\left(\begin{array}{cc} 0 &
    {\rm sign}(q)\sqrt{|q|} \\ -\frac{1}{\sqrt{|q|}} & 0\end{array}\right) .
\]

The irreducible unitary corepresentations have been determined in
\cite{woronowicz87,woronowicz87b,vaksman+soibelman88,masuda+al88,koornwinder89}.
For each non-negative half-integer
$s\in\frac{1}{2}\mathbb{Z}_+$ there exists a $2s+1$-dimensional irreducible
unitary corepresentation $u^{(s)}=(u^{(s)}_{k\ell})_{1\le k,\ell\le 2s+1}$ of
$SU_q(2)$, which is unique up to unitary equivalence and contragredient to
itself. The map $\phi\to\phi^\dag$ on the dual discrete algebraic quantum
group maps the summands $\hat{V}_{\{s\}}\cong M_{2s+1}(\mathbb{C})$ in the
decomposition \eqref{eq-decomp} to themselves and takes the form
\[
A^\dag = Q\overline{A}Q^{-1}
\]
where $A\mapsto\overline{A}$ denotes entry-wise complex conjugation and
\[
Q= \left((-1)^jq^{j-1}\delta_{j,n-k+1}\right)_{1\le j,k\le 2s+1}\in
  M_{2s+1}(\mathbb{C}).
\]
For example  for $s=\frac{1}{2}$, the fundamental corepresentation $u^{(1/2)}=u$, $Q=F_q$, and
\[
R_{1/2} \cong \left\{\left(\begin{array}{cc} a & -q\overline{b} \\ b & \overline{a}\end{array}\right);a,b\in\mathbb{C}\right\}.
\]
For $q>0$, $\mathbb{H}\ni \alpha+\beta I+\gamma J+\delta K \mapsto
\left(\begin{array}{cc} \alpha+i\beta & \sqrt{q}(-\gamma+i\delta) \\
    \frac{1}{\sqrt{q}}(\gamma+i\delta) &  \alpha-i\beta\end{array}\right)\in
R_{1/2}$ defines an isomorphism of real algebras and $R_{1/2}$ contains no
non-zero nilpotent elements.

For $q<0$, $R_{1/2}$ is isomorphic to $M_2(\mathbb{R})$ and contains nilpotent
elements, e.g., $\left(\begin{array}{cc} \sqrt{q} & -q \\ 1 &
    -\sqrt{q} \end{array}\right)$.

The higher dimensional irreducible unitary corepresentations always give
non-zero nilpotent hermitian functionals.
\end{example}

\section{Kac property and the `local' structure of quantum groups in the quantum $DS$-family}

As an application of Theorem \ref{thm-structure} we will now show that  compact quantum groups in the quantum
$DS$-family  are necessarily \emph{of Kac type}, i.e.\ its Haar state is a trace. Recall that the Haar state of a
compact quantum group $(\mathsf{A},\Delta)$ is tracial if and only if the antipode $S$ on $\mathcal{A}$ is
involutive, see \cite[Theorem 1.5]{woronowicz98}. This is the case if and only
if for any unitary corepresentation $u=(u_{jk})_{1\le j,k\le n}$ the
corepresentation $\bar{u}= (u_{jk}^*)_{1\le j,k\le n}$ obtained by taking
adjoints component-wise is again unitary.

Hence being of Kac type is in a sense a `local' property, which
will be clear from the proof of Theorem \ref{thm-kac-type}. We first need to recall a few more facts and
definitions.

%Group $C^*$-algebras of discrete groups and the algebra of continuous functions on a compact group are examples of compact quantum groups of Kac type.

 A compact quantum group $(\mathsf{A},\Delta)$ is called a \emph{compact matrix quantum group} if it has a
finite-dimensional corepresentation $u$ whose coefficients generate $\mathsf{A}$ as a C$^*$-algebra. It follows
from \cite[Proposition 3.7]{maes+vandaele98} that the C$^*$-algebra $\mathsf{A}(u)={\rm C}^*\left(\{u_{jk}:1\le
j,k\le n\}\right)$ generated by the coefficients of any unitary corepresentation $u\in M_n(\mathbb{C})$ of a
compact quantum group $(\mathsf{A},\Delta)$ is a compact matrix quantum group with the restriction of the
coproduct of $\mathsf{A}$. We will call $(\mathsf{A}(u),\Delta|_{\mathsf{A}(u)})$ the \emph{quotient quantum
group} of $(\mathsf{A},\Delta)$ induced by $u$. Equivalent corepresentations clearly induce isomorphic quotient
quantum groups.

\begin{theorem}\label{thm-kac-type}
Let $(\mathsf{A},\Delta)$ be a compact quantum group in the quantum $DS$-family. Then $(\mathsf{A},\Delta)$ is of Kac type.
\end{theorem}
\begin{proof}
Assume that $(\mathsf{A},\Delta)$ is in the quantum $DS$-family. It is sufficient to show that the square of the
antipode acts identically on the coefficients of the irreducible unitary corepresentations of
$(\mathsf{A},\Delta)$. By Theorem \ref{thm-structure}, the irreducible unitary corepresentations of
$(\mathsf{A},\Delta)$ are have dimension one or two.

Let us first consider the one-dimensional corepresentations. If $a$ is
the coefficient of a one-dimensional unitary corepresentation of
$(\mathsf{A},\Delta)$, then $a$ is group-like, i.e.\ $\varepsilon(a)=1$ and
$\Delta(a)=a\otimes a$. Therefore $S(a)=a^*=a^{-1}$ and $S^2(a)=a$.

Let now $u=(u_{jk})_{1\le j,k\le 2}\in M_2(\mathsf{A})$ be a two-dimensional irreducible unitary corepresentation
of $(\mathsf{A},\Delta)$. We will show that the quotient quantum group $(\mathsf{A}(u),\Delta|_{\mathsf{A}(u)})$
is a quantum subgroup of $SU_q(2)$ for some $0<q\le 1$.

By Theorem \ref{thm-structure} and Proposition \ref{prop-char}, there exists
an invertible matrix $Q\in M_2(\mathbb{C})$ such that
\[
\overline{u} = TuT^{-1},
\]
with $T=\bar{Q}\left(\begin{array}{cc} 0 & 1 \\ -1 & 0\end{array}\right)
Q^{-1}$. The matrix $T$ satisfies the relation $T\overline{T}=-I$, so by
\cite[Equation (5.4)]{bichon+rijdt+vaes06}, there exist $0<q\le 1$ and a
unitary matrix $U\in M_2(\mathbb{C})$ such that
\[
T= U^t \left(\begin{array}{cc} 0 & q \\ -\frac{1}{q} & 0 \end{array}\right) U.
\]
Let $v=UuU^*$, then clearly $v$ is a two-dimensional irreducible unitary corepresentation of
$(\mathsf{A},\Delta)$ and $(\mathsf{A}(v),\Delta|_{\mathsf{A}(v)})=(\mathsf{A}(u),\Delta|_{\mathsf{A}(u)})$.
Furthermore, $v$ satisfies the relation
\begin{eqnarray*}
\overline{v} &=& \overline{UuU^*}=\overline{U}\overline{u}U^t
=\overline{U}TuT^{-1}U^t \\
&=&  \overline{U}U^t\left(\begin{array}{cc} 0 & q \\
    -\frac{1}{q} & 0 \end{array}\right)UuU^{-1} \left(\begin{array}{cc} 0 & q \\
    -\frac{1}{q} & 0 \end{array}\right)^{-1}(U^t)^{-1}U^t \\
&=& \left(\begin{array}{cc} 0 & q \\
    -\frac{1}{q} & 0 \end{array}\right)v \left(\begin{array}{cc} 0 & q \\
    -\frac{1}{q} & 0 \end{array}\right)^{-1},
\end{eqnarray*}
i.e.\ the coefficients of $v$ satisfy the defining relations of the universal
orthogonal quantum group $A_0\left(\begin{array}{cc} 0 & q \\
    -\frac{1}{q} & 0 \end{array}\right)\cong SU_{q^2}(2)$, cf.\
\cite{vandaele+wang96} or Example \ref{exa-suq}. If we denote by $w=(w_{jk})_{1\le j,k\le 2}$ the generators of
$SU_{q^2}(2)$, then $w_{jk}\mapsto v_{jk}$ defines a surjective morphism of compact quantum groups from $SU_q(2)$
to $(\mathsf{A}(v),\Delta|_{\mathsf{A}(v)})$, i.e.\ $(\mathsf{A}(v),\Delta|_{\mathsf{A}(v)})$ is a quantum
subgroup of $SU_{q^2}(2)$.

In Example \ref{exa-suq} we have seen that the Haar state on $SU_{q^2}(2)$ admits a non-trivial square root, so
$(\mathsf{A}(u),\Delta|_{\mathsf{A}(u)})$ has to be a proper quantum subgroup of $SU_{q^2}(2)$. For $q^2\not=1$,
the quantum subgroups of $SU_{q^2}(2)$ are the torus $\mathbb{T}$ and its subgroups, cf.\ \cite{podles95} or also
\cite{franz+skalski+tomatsu09}. These are classical groups and therefore of Kac type. For $q^2=1$ we get
$SU_1(2)$, which is a also classical group and of Kac type.

It follows that the quotient quantum group $(\mathsf{A}(u),\Delta|_{\mathsf{A}(u)})$ induced by any irreducible
unitary corepresentation $u$ is of Kac type, therefore we have $S^2={\rm id}$ on the dense $*$-Hopf algebra
contained in $\mathsf{A}$, and $(\mathsf{A},\Delta)$ is of Kac type.
\end{proof}

The proof of the above theorem shows that the structure of quotient quantum groups induced by two-dimensional
irreducible corepresentations is in fact quite rigid. This is formalised in the following corollary.

\begin{corollary}\label{cor-kac1}
A compact quantum group $(\mathsf{A},\Delta)$ belongs to the quantum $DS$-family if and only if the following two conditions are satisfied:
\begin{rlist}
\item $(\mathsf{A},\Delta)$ admits only one- and two-dimensional irreducible unitary corepresentations;
\item the quotient quantum groups
$(\mathsf{A}(u),\Delta|_{\mathsf{A}(u)})$ of $(\mathsf{A},\Delta)$ induced by its two-dimensional irreducible
unitary corepresentations are isomorphic to $C(H)$, where $H$ is the eight-element group of unit quaternions.
\end{rlist}

\end{corollary}
\begin{proof}

Let $(\mathsf{A},\Delta)$ belong to the quantum $DS$-family. Condition (i) follows from  Theorem \ref{thm-structure} and Proposition \ref{prop-char}. Let then $u$ be a two-dimensional irreducible unitary corepresentation of $(\mathsf{A},\Delta)$. The second part of the proof of Theorem \ref{thm-kac-type} implies that $(\mathsf{A}(u),\Delta|_{\mathsf{A}(u)})$
is isomorphic to a proper quantum subgroup of $C(SU_q(2))$ for some $q\in(0,1]$, so, due to \cite{podles95}, also
to a proper quantum subgroup of $C(SU(2))$.

Condition (ii) then follows by inspection of the subgroups of $SU(2)$ (see, e.g., \cite{podles95}), since $H$ is the
only subgroup of $SU(2)$ which is in the $DS$-family and which  has a two-dimensional irreducible unitary
representation.

Conversely, if $(\mathsf{A},\Delta)$ satisfies (i)-(ii) above, then each of  its irreducible unitary corepresentations verifies one of the conditions in Proposition \ref{prop-char}. Thus Theorem \ref{thm-structure} ends the proof.

\end{proof}

The above corollary can be interpreted as describing the `local' structure of the elements of the quantum
$DS$-family -- we know that they only admit one- and two-dimensional irreducible corepresentations and have now
the full understanding of what types of quotient quantum groups are generated by each of the irreducible
corepresentations (for one-dimensional corepresentations the resulting quotients are just the algebras of
functions on cyclic groups). It has one other important consequence.
Recall the Nichols-Zoeller theorem that states that the dimension of a finite-dimensional Hopf algebra is divisible by the dimensions of its Hopf subalgebras, cf.\ \cite{nichols+zoeller89}. Together with the above corollary it implies the following result. %We have just shown that any non-cocommutative compact quantum group whose Haar state doesn not admit a non-trivial square root contains a Hopf subalgebra isomorphic to $C(H)$. This implies the following condition on the dimensions of finite-dimensional non-cocommutative compact quantum group whose Haar state doesn not admit a non-trivial square root.

\begin{corollary} \label{cor-NZ}
Let $(\mathsf{A},\Delta)$ be a non-cocommutative compact quantum group in the quantum $DS$-family. If
$\mathsf{A}$ is finite-dimensional, then its dimension is divisible by eight.
\end{corollary}

We finish this section by describing the structure of the algebra of functions on $H$ in greater detail; this
will be of use in the next section.

%In this section we continue to investigate the structure of the quantum $DS$-family. We already know that it contains only compact quantum groups of Kac type. It turns out that  the eight-element group of unit quaternions again plays a central role, like in Diaconis and Shahshahani's
% classification. We shall need the following notations.

\begin{example}\label{exa-quaternion}
Denote by $\pm 1,\pm I, \pm J,\pm K$ the eight unit quaternions, with the relations
$I^2=J^2=K^2=-1$, $I\cdot J= K$, $J\cdot I=-K$, etc. Denote by $\lambda_g$ and $\mathbf{1}_{\{g\}}$, $g\in \{\pm 1,\pm I, \pm
J,\pm K\}$ the corresponding bases for $C^*(H)\cong \mathbb{C}H$ and $C(H)$. Besides the constant function $\mathbf{1}_H$, $H$ has
three more one-dimensional irreducible unitary representations, which are uniquely determined by
\begin{gather*}
\sigma_I(I)=1, \quad \sigma_I(J)= -1, \\
\sigma_J(I)=-1, \quad \sigma_J(J)= 1, \\
\sigma_K(I)=-1, \quad \sigma_K(J)= -1.
\end{gather*}
Furthermore, $H$ has, up to unitary equivalence, a unique two-dimensional irreducible unitary representation
$\pi:H\to M_2(\mathbb{C})$ of $H$ (or, equivalently, corepresentation of $C(H)$), given by
\[
\pi(I)= \left(\begin{array}{cc}
0 & 1 \\
-1 & 0
\end{array}\right) \quad\mbox{ and }\quad \pi(J)= \left(\begin{array}{cc}
0 & i \\
i & 0
\end{array}\right).
\]
Denote by $\pi_{jk}\in C(H)$, $1\le j,k\le 2$ the matrix elements of $H$ w.r.t.\ to the standard basis of $\mathbb{C}^2$. Then $\{\mathbf{1}_{H},\sigma_I,\sigma_J,\sigma_K,\pi_{11},\pi_{12},\pi_{21},\pi_{22}\}$ is a basis of $C(H)$. We set
\[
C(H)_0 = {\rm span}\{\mathbf{1}_{H},\sigma_I,\sigma_J,\sigma_K\}
\quad\mbox{ and }\quad
C(H)_1 = {\rm span}\{\pi_{11},\pi_{12},\pi_{21},\pi_{22}\}.
\]
We have $\Delta \pi_{jk} = \sum_{\ell=1}^2 \pi_{j\ell}\otimes\pi_{\ell k}$ and $\Delta g = g\otimes g$ for the one-dimensional unitary corepresentation of $C(H)$. Furthermore, on can check that the tensor product of the two-dimensional representation of $H$ decomposes into a direct sum of the four one-dimensional representations. From these observations follows that the decomposition $C(H)=C(H)_0\oplus C(H)_1$ defines a $\mathbb{Z}_2$-grading of $C(H)$, i.e.\ we have
\begin{eqnarray*}
\big(C(H)_j\big)^* &\subseteq& C(H)_j, \\
C(H)_j\cdot C(H)_k&\subseteq& C(H)_{j+k}, \\
\Delta C(H)_j &\subseteq& C(H)_j\otimes C(H)_j
\end{eqnarray*}
for $j,k\in\mathbb{Z}_2$. Equivalently, the map $d:C(H)\to C(H)\otimes \mathbb{C}\mathbb{Z}_2$ defined by
\[
d(u) = u_0\otimes \delta_0 + u_1\otimes \delta_1
\]
for $u=u_0+u_1$ with $u_0\in\mathsf{A}_0$, $u_1\in\mathsf{A}_1$, and $\delta_g$, $g\in \mathbb{Z}_2$ the standard
basis of $C^*(\mathbb{Z}_2)\cong\mathbb{C}\mathbb{Z}_2$, defines a coaction of $C^*(\mathbb{Z}_2)$ on $C(H)$.

In yet another way we can interpret the grading described above as resulting from the fact that $\mathbb{C}\mathbb{Z}_2$ is a quotient Hopf $^*$-algebra of $C(H)$; this will be used in Example \ref{exam-bicrossed}.
\end{example}

\section{Combining the `local' structure of the compact quantum groups in the quantum $DS$-family into the global one and genuinely quantum examples}

In the last section we showed that if $(\mathsf{A},  \Delta)$ is in the quantum $DS$-family, then we can
completely determine the quotient quantum subgroups induced by individual irreducible corepresentations of
$\mathsf{A}$. Here we show that they can be combined in a non-trivial way to provide the examples which belong to
the quantum $DS$-family and are neither commutative nor cocommutative.

It is easy to see that the quotient quantum subgroup induced by an arbitrary number of one-dimensional
irreducible corepresentations of a compact quantum group is always cocommutative, so of the form $C^*(\Gamma)$
for some (discrete) group $\Gamma$. We therefore begin our analysis by combining a one-dimensional
corepresentation with a two-dimensional one.

\begin{proposition}\label{prop-kac2}
Let $(\mathsf{A},\Delta)$ be in the quantum $DS$-family. Let $u\in
M_2(\mathsf{A})$ be a two-dimensional irreducible unitary corepresentation and $g\in\mathsf{A}$ a one-dimensional unitary
corepresentation. Then there exists a unitary matrix $U\in M_2(\mathbb{C})$ s.t.\
\[
\left(\begin{array}{cc} g u_{11} & gu_{12} \\ g u_{21} & g
    u_{22}\end{array}\right) = U \left(\begin{array}{cc}
    u_{11}g^{-1} & u_{12}g^{-1} \\ u_{21}g^{-1} & u_{22}g^{-1}\end{array}\right)U^*.
\]

\end{proposition}
\begin{proof}
Multiplying $u$ by $g$, we get a  two-dimensional irreducible unitary
corepresentation $gu=\left(\begin{array}{cc} g u_{11} & gu_{12} \\ g u_{21} & g
    u_{22}\end{array}\right)$. By Proposition \ref{prop-char}, $u$ and $ug$
are contragradient to themselves. Since $(\mathsf{A},\Delta)$ is Kac,
$\bar{u}$ and $ \overline{gu}=\bar{u}g^{-1}$ are unitary. Therefore the pairs
$u$ and $\bar{u}$, and $gu$ and $\bar{u}g^{-1}$, are unitarily equivalent,
which implies that $gu$ and $ug^{-1}$ are also unitarily equivalent.
\end{proof}

Let $u$ be a two-dimensional irreducible unitary corepresentation of a compact quantum group $(\mathsf{A},\Delta)$
in the quantum $DS$-family, and let $g_1,\ldots,g_n$ be one-dimen\-sional unitary corepresentations that
do not belong to $\mathsf{A}(u)$. The above proposition  suggests that $(\mathsf{A}(u),\Delta|_{\mathsf{A}(u)})$ and
$(\mathsf{A}(g_1\oplus\cdots\oplus g_n),\Delta|_{\mathsf{A}(g_1\oplus\cdots\oplus g_n)})$ are a matched pair, and that quotient
quantum group $(\mathsf{A}(v),\Delta|_{\mathsf{A}(v)})$ generated by the direct sum $v=u\oplus g_1\oplus \cdots \oplus g_n$ is
given by a bicrossproduct of $(\mathsf{A}(u),\Delta|_{\mathsf{A}(u)})$ and $(\mathsf{A}(g_1\oplus\cdots\oplus
g_n),\Delta|_{\mathsf{A}(g_1\oplus\cdots\oplus g_n)})$. We use this idea to construct examples of noncommutative,
noncocommutative compact quantum groups of Kac type whose Haar state admits no non-trivial square roots.

\begin{example} \label{exam-bicrossed}
Let $\Gamma$ be a commutative discrete group and $C(H)=C(H)_0\oplus C(H)_1$ the algebra of functions on the eight-element group
of unit quaternions, with the $\mathbb{Z}_2$-grading introduced in Example \ref{exa-quaternion}. Denote the quotient Hopf $^*$-algebra morphism from $C(H)$ onto
$\mathbb{C}\mathbb{Z}_2$ by $p$. Note that $\mathbb{Z}_2$ acts on $\mathbb{C}\Gamma$ via the standard (period 2) automorphism, mapping $\gamma$ to $\gamma^{-1}$ and let $\tilde{\alpha}:  \mathbb{C}\mathbb{Z}_2 \otimes \mathbb{C}\Gamma\to \mathbb{C}\Gamma$ denote the corresponding map on the level of Hopf $^*$-algebras,
so that
\[ \tilde{\alpha} (\lambda_0 \otimes \gamma) = \lambda_0 \otimes \gamma, \;\;  \tilde{\alpha} (\lambda_1 \otimes \gamma) = \lambda_1 \otimes \gamma^{-1}.\]

Define  now the action $\alpha:C(H)\otimes \mathbb{C}\Gamma\to \mathbb{C}\Gamma$ by the formula
\[\alpha = \tilde{\alpha} \circ (p \otimes \textup{id}_{\mathbb{C}\Gamma}).\]
It turns $\mathbb{C}\Gamma$ into a left $C(H)$-module algebra and comodule algebra, as follows from the fact that  $\tilde{\alpha}$ has this property and $p$ is a Hopf algebra morphism.

Consider now the crossed (smashed) product construction of $\mathbb{C}\Gamma \rtimes C(H)$, following for example \cite[Proposition 1.6.6]{Majidbook}.
The vector space
$K=\mathbb{C}\Gamma \odot C(H)$ can be turned into an algebra with the multiplication (we use here and below the Sweedler's notation for coproducts)
\[
(a\otimes u)(b\otimes v) = a\alpha(u_{(1)} \otimes b) \otimes  u_{(2)} v
\]
($u,v\in C(H)$, $a,b\in \mathbb{C}\Gamma$)
Moreover the proof of \cite[Proposition 6.2.1]{Majidbook} implies that with the standard tensor coproduct
\[
\Delta_K (a\otimes u) = a_{(1)}\otimes u_{(1)} \otimes a_{(2)}\otimes u_{(2)},
\]
 $K$ becomes a Hopf algebra -- note that although $C(H)$ is not cocommutative, so formally we cannot apply directly
\cite[Proposition 6.2.1]{Majidbook}, in fact the action $\alpha$ is defined so that it only `sees' the cocommutative quotient $\mathbb{C}\mathbb{Z}_2$. This explains why we still obtain the desired result.

Using further  the definition of $\alpha$ we get
\begin{eqnarray*}
(a\otimes u_j)(b\otimes v) =   a \varepsilon(u_{j(1)}) S^j(b) \otimes u_{j(2)} v  = a S^j(b) \otimes u_j v
\end{eqnarray*}
for $j\in\{0,1\}$, $u_j\in C(H)_j$, $v\in C(H)$, $a,b\in\mathbb{C}\Gamma$. One can check that $K$ becomes further a $*$-Hopf algebra with
\[
( a\otimes u_j)^* = \big((a\otimes 1)(1\otimes u_j)\big)^* = (1\otimes u_j^*)(a^*\otimes 1) =  S^j(a^*) \otimes u_j^*
\]
for $j\in\{0,1\}$, $u_j\in C(H)_j$, $a\in\mathbb{C}\Gamma$. Since Haar states are invariant under the respective antipodes, we
see the tensor product of the Haar states on $C^*(\Gamma)$  and $C(H)$ defines a normalized positive integral on $K$, i.e.\ $K$ is
an algebraic compact quantum group in the sense of \cite{vandaele98,vandaele03}.

If $u\in M_n\big(C(H)\big))$ is an irreducible unitary corepresentation of $C(H)$ and $g\in \Gamma$, then $u\otimes g$ is unitary
in $K$, e.g., in the simplified notation,
\[
(g\otimes u) (g\otimes u)^* = (g\otimes u)( S(g^*)\otimes u^*) = (g\otimes u)(g\otimes u^*) =  gS(g) \otimes uu^*=
1\otimes 1.
\]
%for $u\in M_2\big(C(H)\big)$ the two-dimensional irreducible unitary corepresentation of $C(H)$.

Therefore, and since the coproduct in $K$ is simply the tensor product of the coproducts of  $\mathbb{C}\Gamma$ and $C(H)$, we see
that the irreducible unitary corepresentations of $K$ are of the form $g\otimes u$, with $u$ an irreducible unitary
corepresentation of $C(H)$, and $g\in \Gamma$. Since their coefficients span $K$, we can deduce that $K$ is the $*$-Hopf algebra
of a compact quantum group, for which we can choose the universal $C^*$-algebra of the $*$-algebra $K$. Let us denote this
quantum group by $ C^*\Gamma\rtimes_\alpha C(H) $.

The construction above can also be viewed as a special case of the double crossed product construction from \cite[Proposition 3.12]{majid90b} or \cite[Example 6.2.12]{Majidbook}.

It is clear that $ C^*\Gamma\rtimes_\alpha C(H) $ has only one- and two-dimensional irreducible unitary corepresentations, and
straight-forward to show that the two-dimensional irreducible corepresentations are self-contragredient and satisfy the condition
in Proposition \ref{prop-char}. Hence the Haar state of $ C^*\Gamma\rtimes_\alpha C(H) $ does not admit any non-trivial square
roots.
\end{example}

\begin{proposition} \label{prop-quantexam}
Let $\Gamma$ be an abelian discrete group which contains elements that are not of order two. The crossed
product $ C^*\Gamma\rtimes_\alpha C(H) $ constructed in Example \ref{exam-bicrossed} is a noncommutative,
noncocommutative compact quantum group in the quantum $DS$-family.
\end{proposition}

Proposition \ref{prop-kac2} describes the interaction of a two-dimensional irreducible corepresentation, say $u$,
of a quantum group in the quantum $DS$-family with one-dimensional ones. It is reflected by certain equivalences
between $u \otimes g$ and $g^{-1} \otimes u$.  In the last part of the paper we discuss certain aspects of the
interaction between different two-dimensional representations. To this end we need to introduce a certain
equivalence relation on the equivalence classes of irreducible corepresentations of a fixed compact quantum
group. Introduce first the notation: if  $(\mathsf{A},\Delta)$ is a compact quantum group, let
$\textup{Irr}(\mathsf{A})$ denote the set of the equivalence classes of irreducible corepresentations of
$\mathsf{A}$ and let $\Gamma_{\mathsf{A}} \subset \textup{Irr}(\mathsf{A})$ denote the equivalence classes of
one-dimensional corepresentations (in other words, group-like elements of $\mathsf{A}$). It is well-known (and
has been used above) that the tensor product of corepresentations provides $\Gamma_{\mathsf{A}}$ with the
structure of a discrete group.

\begin{proposition}
Let $(\mathsf{A},\Delta)$ be a compact quantum group. The relation $\approx_{\,\Gamma}$ on
$\textup{Irr}(\mathsf{A})$ given by the formula
\[ u \approx_{\,\Gamma} v \;\;\; \textup { if } \;\; \exists_{\gamma \in \Gamma_{\mathsf{A}}} \, u = v \otimes \gamma\]
is an equivalence relation.
\end{proposition}
\begin{proof}
Easy check, essentially a consequence of the fact that $\Gamma_{\mathsf{A}}$ forms a group and associativity of
the tensor operation for (not necessarily irreducible) corepresentations of $\mathsf{A}$.
\end{proof}

The set of equivalence classes in $\textup{Irr}(\mathsf{A})$ with respect to the relation $\approx_{\,\Gamma}$
will be denoted $\textup{Irr}(\mathsf{A})/\approx_{\,\Gamma}$ and for $u\in \textup{Irr}(\mathsf{A})$ the
corresponding equivalence class $\textup{Irr}(\mathsf{A})/\approx_{\,\Gamma}$ will be denoted $[u]_{\approx}$.
Note that all one-dimensional corepresentations form an equivalence class with respect to the relation
$\approx_{\,\Gamma}$, to be denoted $[1]_{\approx}$.

\begin{theorem}
Let $(\mathsf{A},\Delta)$ be in the quantum $DS$-family. Then the set $\textup{Irr}(\mathsf{A})/\approx_{\Gamma}$
is equipped with a well-defined product, which is given by the following condition: for $u,v,w \in
\textup{Irr}(\mathsf{A})$
\begin{equation} \label{prodirr} [u]_{\approx} \cdot [v]_{\approx} = [w]_{\approx}
\;\;\; \textup { if } \;\; \exists_{\gamma \in \Gamma_{\mathsf{A}}} \; u \otimes v \succeq w \otimes
\gamma.\end{equation}
 Moreover the pair $(\textup{Irr}(\mathsf{A})/\approx_{\,\Gamma}, \cdot)$ forms an abelian
group, in which each non-trivial element has order 2.
\end{theorem}
\begin{proof}
We need to check first that the product in the formula \eqref{prodirr} is well defined. If $\gamma, \gamma' \in
\textup{Irr}(\mathsf{A})$ are one-dimensional, then so is $\gamma\otimes \gamma'$ and according to the notation
of \eqref{prodirr} and that introduced before the statement of the theorem we have $[1]_{\approx} \cdot [1]_{\approx} = [1]_{\approx}$. If $u\in \textup{Irr}(\mathsf{A})$ is two-dimensional and
$\gamma\in \Gamma_{\mathsf{A}}$, then $u \otimes \gamma$ is irreducible and equivalent to $u$ with respect to
$\approx_{\,\Gamma}$, so $[u]_{\approx} \cdot [1]_{\approx} = [u]_{\approx}$. Similarly $\gamma \otimes u$ is
irreducible. Moreover, as both  $u$ and $\gamma \otimes u$ are self-contragredient by Proposition
\ref{prop-char}, $\gamma \otimes u = u \otimes \gamma^{-1}$ (recall that the equality here is understood in terms
of the usual equivalence classes of irreducible representations).

The only non-trivial case is that of $u,v \in \textup{Irr}(\mathsf{A})$  both two-dimensional. Observe that due
to Corollary \ref{cor-kac1} and the discussion in Example \ref{exa-quaternion} the tensor product $u \otimes u$ is a
four-dimensional corepresentation decomposing into 4 one-dimensional corepresentations, including the trivial
one.  We will distinguish two-possibilities: first assume that $u \approx_{\,\Gamma} v$. Then there is some
$\gamma \in \Gamma_{\mathsf{A}}$ such that $u =v \otimes \gamma= \gamma^{-1} \otimes v$. Hence $u \otimes v = \gamma^{-1}\otimes
(v \otimes v)$ is a direct sum of four one-dimensional corepresentations (trivially equivalent to each other with
respect to  $\approx_{\,\Gamma}$). This can be rephrased by writing $[u]_{\approx} \cdot [u]_{\approx} =
[1]_{\approx}$. It remains to consider the possibility of $u \not\approx_{\,\Gamma} v$.  Then $u\otimes v$ is a
four-dimensional, necessarily reducible corepresentation. Suppose that $u \otimes v$ contains a one-dimensional
corepresentation, say $\gamma$. Then $(\gamma^{-1} \otimes u) \otimes v$ contains a trivial corepresentation, and
as both $\gamma^{-1} \otimes u$ and $v$ are irreducible, this would mean that $\gamma^{-1} \otimes u=v^c=v$, so
$u \approx_{\,\Gamma} v $ -- contradiction. Thus $u \otimes v$ decomposes into a direct sum of two
two-dimensional irreducible corepresentations, say $w$ and $w'$. To assure that the product in \eqref{prodirr} is
well defined we need to prove that $w \approx_{\,\Gamma} w'$ (strictly speaking we also need to show that
$[w]_{\gamma}$ depends only on the $\approx_{\,\Gamma}$-equivalence classes of $u$ and $v$, but this is easy to
see). Tensor the corepresentation $u \otimes v$ on the left with $u$. Then it decomposes into four
two-dimensional irreducible corepresentations, each $\approx_{\,\Gamma}$-equivalent to $v$. Hence in particular
$u\otimes w$ as a subrepresentation of $u \otimes v$  decomposes into two-dimensional irreducible corepresentations $\approx_{\,\Gamma}$-equivalent to
$v$, say $\gamma_1 \otimes v,$ and $\gamma_2 \otimes v$. Tensor the formula $u\otimes w=(\gamma_1 \otimes v)
\oplus (\gamma_2\otimes v)$ again on the left with $u$. Then on the left we obtain the direct sum of four
two-dimensional corepresentations, each $\approx_{\,\Gamma}$-equivalent to $w$, and on the right the direct sum
of four two-dimensional corepresentations, two of which are $\approx_{\,\Gamma}$-equivalent to $w$, and two are
$\approx_{\,\Gamma}$-equivalent to $w'$. Hence $w$ and $w'$ are $\approx_{\,\Gamma}$-equivalent and the proof of
the main part of the theorem is finished.

As to the fact that the product $\cdot$ gives $\textup{Irr}(\mathsf{A})/\approx_{\Gamma}$ the group structure
described in the theorem, it suffices to observe that $\cdot$ inherits associativity from the usual associativity
of tensor products of corepresentations and that the first part of the proof shows that $[1]_{\approx}$ is the
neutral element for $\cdot$ and each element in $\textup{Irr}(\mathsf{A})/\approx_{\Gamma}$ is its own inverse.
\end{proof}

The above theorem implies that if $(\mathsf{A},\Delta)$ is in the quantum $DS$-family then $(\textup{Irr}(\mathsf{A})/\approx_{\,\Gamma}, \cdot)$ is a direct sum of (possibly infinitely many) copies of $\mathbb{Z}_2$. For cocommutative $(\mathsf{A},\Delta)$ the group $(\textup{Irr}(\mathsf{A})/\approx_{\,\Gamma}, \cdot)$ is trivial. For a compact group $G$ in the $DS$-family the group $(\textup{Irr}(C(G))/\approx_{\,\Gamma}, \cdot)$ is either trivial or a two-element group, depending on whether $G$ contains the group of unit quaternions. If $(\mathsf{A},\Delta)$ is the compact quantum group constructed in Example \ref{exam-bicrossed}, then again the group $(\textup{Irr}(\mathsf{A})/\approx_{\,\Gamma}, \cdot)$ is isomorphic to $\mathbb{Z}_2$. It is therefore natural to seek the answer to the following open question.

\begin{problem}\label{prob-bigIrrGamma}
Does there exist $(\mathsf{A},\Delta)$ in the quantum $DS$-family such that the group  $ (\textup{Irr}(\mathsf{A})/\approx_{\,\Gamma}, \cdot)$ has more than two elements?
\end{problem}

It is natural to seek for such an object among finite-dimensional Kac
algebras, exploiting the existing classification of low-dimensional examples
(see \cite{Izumi+Kosaki2002}). Corollary \ref{cor-NZ} implies that the dimension of such a Kac algebra would have to be divisible by $8$. To allow for two distinct classes of two-dimensional irreducible corepresentations we need the dimension to be at least $16$. The case by case analysis of the form of the Grothendieck rings of 16-dimensional Kac algebras listed in \cite{kashina00} implies that none of these algebras can provide a positive answer to the question asked in Problem \ref{prob-bigIrrGamma}. Hence the lowest dimension for the Kac algebra that would answer the question in Problem \ref{prob-bigIrrGamma} is equal 24.

\section*{Acknowledgements}
Part of this work was carried out while U.F.\ was visiting the Graduate School
of Information Sciences of Tohoku University as Marie-Curie fellow. He would
like to thank Professors Nobuaki Obata, Fumio Hiai, and the other members of
the GSIS for their hospitality.

%\bibliographystyle{alpha}
%\bibliography{../../Bibliography/bibliography}

%\end{document}

\end{document}